\documentclass[9.5pt]{amsart}
\usepackage{amssymb}

\newtheorem{lemma}{Lemma}
\newtheorem{prop}[lemma]{Proposition}
\newtheorem{cor}[lemma]{Corollary}
\newtheorem{thm}[lemma]{Theorem}
\newtheorem{fact}[lemma]{Fact}

\newtheorem*{mainthm}{Main Theorem}
\newtheorem{thm?}[lemma]{Theorem?}

\DeclareFontEncoding{OT2}{}{} 

\title{Bounds for torsion on abelian varieties with integral moduli}
\author{Pete L. Clark}
\address{1126 Burnside Hall \\ Department of Mathematics and Statistics \\
McGill University \\ 805 Sherbrooke West \\ Montreal, QC, Canada H3A 2K6}
\email{clark@math.mcgill.ca}

\newcommand{\F}{\ensuremath{\mathbb F}}

\newcommand{\Q}{\ensuremath{\mathbb Q}}

\newcommand{\Z}{\ensuremath{\mathbb Z}}

\newcommand{\ra}{\ensuremath{\rightarrow}}

\newcommand{\Ker}{\operatorname{Ker}}

\newcommand{\ord}{\operatorname{ord}}
\newcommand{\unr}{\operatorname{unr}}

\newcommand{\PP}{\ensuremath{\mathbb P}}

\newcommand{\Aut}{\operatorname{Aut}}
\newcommand{\tors}{\operatorname{tors}}

\newcommand{\End}{\operatorname{End}}

\newcommand{\OO}{\mathcal{O}}

\newcommand{\Gm}{\mathbb G_m}
\newcommand{\Ga}{\mathbb G_a}

\begin{document}

\maketitle
\noindent
\begin{abstract} 
We give a function $F(d,n,p)$ such that if $K/\Q_p$ is a degree $n$ field extension and
$A/K$ is a $d$-dimensional abelian variety with potentially good reduction, 
then $\#A(K)[\tors] \leq F(d,n,p)$.  Separate attention is given to the prime-to-$p$ torsion
and to the case of purely additive reduction.  These latter bounds 
are applied to classify rational torsion on CM elliptic curves over number fields of
degree at most $3$, on elliptic curves over $\Q$ with integral $j$ (recovering a theorem of Frey), and on 
abelian surfaces 
over $\Q$ with integral moduli.  In the last case, our efforts leave us
with $11$ numbers which may, or may not, arise as the order of the full torsion subgroup.  The largest such number is
$72$.
\end{abstract}
\section{Introduction}
\noindent
Let $K$ be a field.  We say that \emph{torsion is strongly bounded} for abelian varieties of 
dimension
$d$ over $K$ if for all finite field extensions $L/K$ and all $d$-dimensional abelian
varieties $A/L$, there exists an $N = N(d,[L:K])$ such that $\#A(L)[\tors] \leq N$. (If we left out the ``strong,'' it would mean that
$N$ is allowed to depend on $L$ and not just its degree over $K$.)  Note that it
is equivalent to require the order of any torsion \emph{point} on $A(L)$ to be bounded
dependent only on $d$ and $[L:K]$.  Also, if $K'/K$ is a finite field extension, then if torsion
is strongly bounded over $K$, it is automatically strongly bounded over $K'$.  
\\ \\
Example 1: Torsion is strongly bounded in every dimension $d$ for abelian varieties
over $K = \mathbb{F}_p$, since the Weil bounds imply that
for any $d$-dimensional abelian variety $A/\F_{p^n}$ we have (as we shall recall in Section 2)
\begin{equation}
\#A(\mathbb{F}_{p^n}) \leq (1+p^{n/2})^{2d}. 
\end{equation}
``Example'' 2: A celebrated theorem of Merel \cite{Merel} asserts that torsion is strongly bounded when $d=1$ and $K = \Q$.
It has been conjectured by many that torsion is strongly bounded for abelian varieties $A/\Q$ in
every dimension $d$, but this seems much beyond present reach.
\\ \\
It follows from Examples 1 and 2 that torsion is strongly bounded for elliptic curves over any finitely
generated field $K$.  However, the uniform boundedness of torsion among all elliptic curves fails over
most other fields of arithmetic interest.  For instance:
\\ \\
Example 3: Torsion is not bounded for elliptic curves over $\Q_p$: on the
Tate curve $E_{p^n} = \Gm/\langle p^n \rangle$, the point $p$ has exact order $n$.
\\ \\
Note that Tate curves are characterized by having split multiplicative reduction.  
On the other hand, it follows
easily from Example 1 that at least the prime-to-$p$ torsion is strongly bounded among
all abelian varieties $A/\Q_p$ with good reduction.  The main result of the present paper is that
strong boundedness holds for a somewhat wider class of abelian varieties over local fields.
\\ \\
First some terminology and notation: In this paper, a \emph{complete field} means a field $K$ complete with
respect to a discrete valuation $v$, and a \emph{local field} is a complete field with
finite residue field, whose characteristic we denote by $p$.  The absolute ramification
index of a local field is $e = v(p) \leq\infty$.  If $A/K$ is an abelian variety over a local field,
then $A(K)[p^{\infty}]$ denotes the subgroup of $p$-primary torsion and $A(K)[\tors]'$ denotes the subgroup of
torsion of order prime to $p$. 
\\ \\
Recall that an abelian variety over a complete field $K$ is said to have \emph{potentially good reduction}
if there exists a finite field extension $L/K$ such that the base change of $A$ to $L$ is the generic
fibre of an abelian scheme over the valuation ring of $L$.  If $R$ is an (arbitrary) Dedekind domain with
quotient field $K$, we will say that an abelian variety $A/K$ has \emph{integral moduli} 
if for every prime ideal $v$ of $R$, the base change of $A$ to the completion of $K$ at $v$ has potentially
good reduction.\footnote{The ``moduli'' are the coordinates of the point induced by $A$ (together with
some choice of polarization) on a suitable Siegel moduli scheme.  However, in the sequel we will only need
to know that an elliptic curve has integral moduli iff its $j$-invariant lies in $R$.}  Sometimes for brevity
we shall speak of an ``IM abelian variety,'' which means an abelian variety defined over a number field
and having integral moduli.
\\ \\
For $p$ a prime number and $D$ a positive integer, define 
\[m_p(D) := \min_{N \geq 3, \ (N,p) = 1} \ord_p(\#GL_D(\Z/N\Z)), \]
It is well-known (e.g. \cite{ST}) that if $p > 2d+1$, $m_p(2d) = 0$.  
\\ \\
For a positive real number $x$, let $[x]$ be the greatest integer which is at most $x$ and 
$[x]_p$ be the largest power of $p$ which is at most $x$.
\begin{mainthm}
Let $K$ be a local field with residue cardinality $q = p^f$ and absolute ramification index $e$.  Let
$A/K$ be a $d$-dimensional abelian variety with potentially good reduction. \\
a) $\#A(K)[\tors]' \leq [(1+\sqrt{q})^{2d}]$. \\
b) Suppose that $K$ is a $p$-adic field (i.e., $e < \infty$).  Then
\[\#A(K)[\tors] \leq [(1+\sqrt{q})^{2d}] [(1+\sqrt{q})^{2d}]_p p^{fd[\frac{e}{p-1}]+2d m_p(2d)}. \]
c) Suppose $K$ is $p$-adic and $A/K$ has purely additive reduction.  Then the only primes 
$\ell \neq p$ dividing $\#A(K)[\tors]$ are those for which $\ell \leq 2d+1$.  Moreover, there
exists a function $Z(d,e)$ -- i.e., independent of the residue characteristic and the inertial
degree -- such that $A(K)[\tors] \leq Z(d,e)$.
\end{mainthm}
\noindent
In particular, torsion is strongly bounded among $d$-dimensional abelian varieties over $\Q_p$ with
potentially good reduction.
\\ \\
Remark: Notice that if $p > \max\{e+1,2d+1\}$, the exponent of the power of $p$ in part b) is zero.
\\ \\
Remark: Essentially the same result has also been proved -- independently, but with almost identical methods -- by Xavier
Xarles.  The only notable difference is that we give a direct proof of part a) and use this to take care of the
prime-to-$p$ torsion in the general case, whereas Xarles' argument uses results of Lorenzini to give bounds on the prime-to-$p$
part of the component group in terms of the additive rank.  However, in our applications we will either have to 
assume the reduction could be good (the worst-case 
scenario!) and use the bound of part a), or we will be able to argue for a ``not especially nasty'' place of purely additive
reduction and use the stronger bound of part c).  It is worth mentioning that we never apply the bound of part b)
directly.
\\ \\
We will give the proof of the Main Theorem in Section 2, beginning in Section 2.1 with the proof of part a).  
Section 2.2 gives a discussion of certain basic results about N\'eron models, formal groups, and 
component groups, which are then used to prove the rest of the theorem in Section 2.3.  In Section 2.4 we discuss open
problems, comparisons with the elliptic curve case, and mention a considerably more general boundedness theorem due to 
Xarles which was conjectured (or rather, ``suggested'') in an earlier draft.  
\\ \\
One consequence of the Main Theorem is that totally indefinite quaternionic Shimura varieties 
with $\Gamma_1(N)$-level structure will, for sufficiently large $N$, fail to have points rational over
any given $p$-adic field.  Because it is somewhat technical to make precise the moduli problem corresponding
to a  $\Gamma_1(N)$-level structure in the higher-dimensional case, we will content ourselves here with 
the (representative) case of Shimura curves over $\Q$, in which case a careful description 
of the moduli problem may be found in \cite{Buzzard}.
\begin{thm}
Let $B/\Q$ be a nonsplit, indefinite rational quaternion algebra with discriminant $D$, let
$N$ be a positive integer prime to $D$, and let $X^D_1(N)/\Q$ be the Shimura curve with
$\Gamma_1(N)$-level structure.  Let $K$ be a $p$-adic field with residue field $\F_q$, and consider 
the curve $X^D_1(N)/K$.  Then there exists $N_0 = N_0([K:\Q_p])$ (but independent of $D$) such 
that for all $N \geq N_0$, $X^D_1(N)(K) = \emptyset$.
\end{thm}
\noindent
Proof: For $N \geq 4$ and $F/\Q$ any field of characteristic zero, to a point
$P \in X^D_1(N)(F)$ we can associate $A/F$ an abelian surface, $\iota: B \hookrightarrow
\End_F(A) \otimes \Q$ a quaternionic multiplication (QM) structure, and $p \in A(F)$ a point of
exact order $N$.  Now suppose $K$ is a local field and $P = (A,\iota,p) \in X^D_1(N)(K)$.  
A QM surface necessarily has potentially
good reduction.  \\ \indent Indeed, any $d$-dimensional abelian variety admitting as a subring of endomorphisms an 
order $\OO$ in a $2d$-dimensional division algebra has potentially good reduction 
(note that this includes also the CM  case).  We sketch the proof: by a theorem of Grothendieck, after a finite base change there
is no additive part; moreover, the character group, if nontrivial, would be the underyling $\Z$-module of a
(necessarily faithful) representation of $\OO$, but its dimension is at most $d$, whereas any nontrivial
representation of $\OO$ has dimension at least $2d$. \\ \indent Thus the result follows immediately from the Main Theorem.
\\ \\
Evidently the Main Theorem can also abe applied to bound the torsion subgroup of an abelian
variety over a number field with integral moduli.  Indeed:
\begin{cor}
Let $A/K$ be a $d$-dimensional abelian variety over a degree $n$ number field $K/\Q$. \\
a) Assume that $A$ has potentially good reduction at places $v_2$,$v_3$ of $K$ over $2$ and
$3$.  Then $\#A(K)[\tors] \leq [(1+2^{n/2})^{2d}][(1+3^{n/2})^{2d}].$ \\
b) Assume that $A$ has potentially good reduction at a place $v_2$ of $K$ over $2$.  \\ Then
the maximum prime order of a torsion point is $[(1+2^{n/2})]^{2d}]$.
\end{cor}
\noindent
Proof: For part a), let $\iota_2: A(K)[\tors] \ra A(K_{v_2})[\tors]'$ be the composite of
the natural embedding of the global torsion points into the local torsion points with projection 
onto
the prime-to-$2$-torsion, and define $\iota_3: A(K)[\tors] \ra A(K_{v_3})[\tors]'$
similarly.  Then the product map
\[\iota_2 \times \iota_3: A(K)[\tors] \ra A(K_{v_2})[\tors]' \times A(K_{v_3})[\tors]' \]
is evidently an injection, so the result follows from Theorem 1a).  Injecting the odd order
torsion into $A(K_{v_2})[\tors]'$ gives part b).  
\\ \\
Remark: The existence of a strong bound on torsion for IM varieties is due to Alice Silverberg
\cite{Silverberg 1992}, \cite{Silverberg 2001}.  She gives the bound:
\begin{equation}
\#A(K)[\tors] \leq [(1+2^{\#GL_{2d}(\Z/3\Z)n/2})(1+3^{\#GL_{2d}(\Z/4\Z)n/2})]^{2d}. 
\end{equation}
As an example of the improvement that the present results provide, in the case of abelian surfaces over $\Q$,
Corollary 2 gives $\#A(\Q) \leq 1815$, while Silverberg's bound is approximately $4.0262 \times 10^{1275357349}.$
\\ \\
But the bound of Corollary 2a) is still visibly far from the truth: because we are 
double counting the prime-to-$6$ torsion, we will in practice get much better bounds by applying part a) of the Main
Theorem prime by prime and 
collating the results.  The point that we want to emphasize is that, at least when the number field degree is small, this
collation process leads to \emph{short lists}\footnote{Readers will have to judge for themselves whether or not the lists we
give are \emph{too} short, in light of the current academic job market, to sustain the analogy.} of possible orders for the 
torsion subgroup -- i.e., it becomes feasible to give serious individual consideration to each of the elements on the list
as to whether or not they actually arise globally.  
\\ \\
Indeed, in some cases the lists are so short that the temptation to try to decide whether or not they are optimal has proved
to be irresistible, and the remainder of the paper -- which has a different and more computational character -- is devoted to
this investigation.  We prove the following results (of which only the last is completely new):
\begin{thm}
Let $E/\Q$ be an elliptic curve with everywhere potentially good reduction.  Then
$E(\Q)[\tors]$ is one of the following groups: 
\begin{equation}
0, \ \Z/2\Z, \ \Z/3\Z, \ \Z/4\Z, \ \Z/6\Z, \ \Z/2\Z \oplus \Z/2\Z. 
\end{equation}
Conversely, all such groups occur among elliptic curves $E/\Q$ with $j = 0$ and $j= 1728$.
\end{thm}
\noindent
\begin{thm}
Let $K$ be a number field of degree $d$, and $E/K$ an elliptic curve which admits complex multiplication over $\overline{K}$.
\\a) If $d = 2$, then $E(K)[\tors]$ is one of the following groups (all of which occur): \\
$\Z/n\Z$ for $n \in \{1,2,3,4,6,7,10\}$, $\Z/2\Z \times \Z/2n\Z$ for $n \in \{1,2,3,4\}$, $\Z/3\Z \times \Z/3\Z$. \\
b) If $d = 3$, then $E(K)[\tors]$ is one of the following groups (all of which occur): \\
$\Z/n\Z$ for $n \in \{1,2,3,4,6,9,14\}$, $\Z/2\Z \times \Z/2\Z$.
\end{thm}
\begin{thm}
Let $A/\Q$ be an abelian surface with everywhere potentially good reduction.  Then the
order of the torsion subgroup
$\#A(\Q)[\tors]$ is one of
\begin{equation}
1-16, \ 18, \ 19, \ 20, \ 22, \ 24, \ 25, \ 28, \ 30, \ 36, \ 48, \ 60, \ 72. 
\end{equation}
\end{thm}
\noindent
Theorem 3 was etablished for CM elliptic curves over $\Q$ by Loren Olson \cite{Olson} and for all
elliptic curves over $\Q$ with integral $j$-invariant by Gerhard Frey \cite{Frey}.  What is notable
is that the purely local methods of the Main Theorem rule out all possibilities except $\Z/5\Z$, which requires
some global argument.  We give a very short proof that $\Z/5\Z$ cannot occur using part c) of the Main Theorem
and the modularity theorem for elliptic curves. 
\\ \\
As for Theorem 4, it is a special case of work done by H.G. Zimmer and his collaborators: indeed, in a sequence of
papers \cite{Z1}, \cite{Z2}, \cite{Z3} they classify the torsion subgroup of an IM elliptic curve
over number fields of degree $2$ and $3$.  The point here is that their work requires massive
computation.  Of course one expects the CM case to be easier; however, by capitalizing on a theorem of J. L. Parish
\cite{Parish}, it turns out that the most naive possible approach gives an apparently reasonable algorithm
for enumerating the torsion subgroups of CM elliptic curves over number fields of prime degree.  This algorithm
uses the solution of the class number one problem -- i.e., that there are precisely $13$ rational integers
which are $j$-invariants of CM elliptic curves -- but does not require any information about higher class numbers.
To give an idea of the relative difficulties involved, it seems to the author that the amount of calculation needed
to deduce the cubic case of Theorem 4 from the classification results of \cite{Z3} (i.e., to decide which of the
IM elliptic curves have complex multiplication) could equal or exceed the calculations necessary to
establish the theorem.
\\ \\
Section 4 is devoted to the proof of Theorem 5.  The methods are purely local (but use more than the Weil bound (1),
which is surprisingly far from being sharp for higher-dimensional abelian varieties over $\F_2$ and $\F_3$).
We end the paper with an incomplete discussion of the values of (5) that we know are attained: $11$ entries
on the list remain in doubt.
\\ \\
This paper grew out of a chapter of the author's 2003 Harvard thesis \cite{Clark}, whose general topic was local
and global points on Atkin-Lehner quotients of Shimura curves, without and (especially) with $\Gamma_0(N)$ and
$\Gamma_1(N)$ level structure.  The author is indebted to his thesis advisor, Barry Mazur, for mathematical and moral support,
and for general torsion-related inspiration.  In particular versions of the Main Theorem, Theorem 1 and Theorem 5 appear there.  
We remark that the possible torsion subgroups of an abelian surface $A/\Q$ with potential quaternionic
multiplication are significantly more restricted than the general IM case: e.g., the maximum prime order of a torsion point
is at most $7$.  These and other essentially quaternionic results will be treated elsewhere.

\section{The proof of the Main Theorem}
\noindent
In Section 2.1 we prove part a) of the Main Theorem.  Section 2.2 records a theorem of McCallum and some basic
facts concerning N\'eron models and torsion points of abelian varieties over local fields.  These facts must
be very well-known, but we are not aware of suitable references, so for the convenience of the less experienced
reader (and the author) we give complete proofs.  In Section 2.3 we complete the proof of the Main Theorem and
in Section 2.4 we discuss an improvement and some open problems.
\\ \\
Throughout this section we use the following notation: $K$ is a complete field, with discrete valuation 
$v$, valuation ring $R$, maximal ideal $\mathfrak{m}$, and residue field $k$.
\subsection{The proof of part a)}
Let $A/K$ be a $d$-dimensional abelian variety over a 
non-Archimedean local field $K$ of finite residue cardinality $q$, a power of $p$.  When $A/K$ has good reduction,
the bound of the theorem is very well known: indeed, the reduction map $A(K)[\tors]' \ra A(\F_q)[\tors]'$
is an injection, so \[\#A(K)[\tors]' \leq \#A(\F_q) = \#\Ker(1-F) = \det(1-F \ | \ T_{\ell}(A)) = 
\prod_{i=1}^{2d} (1-\lambda_i), \]
where the $\lambda_i$'s are the characteristic roots of Frobenius acting on an $\ell$-adic Tate module.  By
the Riemann hypothesis, $||\lambda_i|| = \sqrt{q}$ for any Archimedean norm $||\ ||$ on $\Q(\lambda_i)$.
So $\#A(\F_q) \leq (1+\sqrt{q})^{2d}$, establishing (1).
\\ \\
The point is that we can get exactly the same bound only assuming that $A/K$ has potential good reduction.
Indeed, according to \cite[p. 498]{ST}, we can find a finite totally ramified extension $L/K$ such that
 $A/L$ has good reduction.  But now the prime-to-$p$ torsion of $A(K)$ is contained in the 
prime-to-$p$ torsion
of $A(L)$, which is isomorphic to the prime-to-$p$ torsion of the good reduction 
$A(\mathbb{F}_q)$, i.e., we have reduced to the case of good reduction without increasing the residue
field, so the bound of $(1+\sqrt{q})^{2d}$ still applies.
\subsection{Some preliminaries} \ \\
We begin by recalling the following ``strong form'' of the criterion of N\'eron-Ogg-Shafarevich.
\begin{thm}(Serre-Tate)
Let $A/K$ be an abelian variety with potentially good reduction over a complete field with residue characteristic 
$p$.  For any $N \geq 3$ and prime to $p$, $A$ acquires good reduction over the $\ell$-torsion field
$L = K(A[N])$.
\end{thm}
\noindent
For $A/K$ a $d$-dimensional abelian variety, we write $\mathcal{A}/R$ for the
N\'eron model of $A/K$.  We will be almost entirely interested in the N\'eron special fibre
$\mathcal{A}_s := \mathcal{A} \otimes_R k$.  
\\ \\
Let $\mathcal{R}: A(K) = \mathcal{A}(R) \ra \mathcal{A}(k) = \mathcal{A}_s(k)$ be the \emph{reduction map} (the first
equality is by the N\'eron mapping property).  The kernel $G^1$ of $\mathcal{R}$ is isomorphic to $\hat{A}(\mathfrak{m})$,
the $K$-analytic Lie group obtained by using the formal group law $\hat{A}$ as a group operation on
$\mathfrak{m}^d$.  The following result bounds the torsion in $G^1$.
\begin{prop}
Let $F(X,Y)$ be a $d$-dimensional formal group law over the valuation ring $R$ of $K$, with associated
``standard'' $K$-analytic Lie group $G^1 = F(\mathfrak{m})$.  Let $H \subset G^1$ be any finite subgroup.  Then: \\
a) The exponent of $H$ divides $p^{[\frac{e}{p-1}]}$.  \\
b) The order of $H$ divides $q^{d[\frac{e}{p-1}]}$.
\end{prop}
\noindent
Proof: We recall three results concerning filtrations on standard groups [LALG, II.IV.9], as simplified by our
assumption that the valuation group is discrete.  Indeed, in this case $G^1$ has a filtration
by normal subgroups $\{G^i = F(\mathfrak{m}^i)\}_{i=1}^{\infty}$.  Let $f_r: x \mapsto x^r$
be the $r$th power map.  It is shown: first, that for all $r$ prime to $p$, $f_r$ is an automorphism
of $G_1$; second that for all $i$, $G^i/G^{i+1} \cong (k,+)^n \cong (\Z/p\Z)^{fn}$; and third, that when 
$i > \frac{e}{p-1}$, $f_p: G^i \ra G^{i+p}$ is an isomorphism.  The result follows from these
facts.
\\ \\
We now recall some facts about the structure of the N\'eron special fibre.  There is a short exact sequence 
\[0 \ra \mathcal{A}_s^0 \ra \mathcal{A}_s \ra \Phi \ra 0, \]
exhibiting the special fibre as an extension of a finite \'etale group scheme $\Phi/k$ (the component
group) by a connected commutative algebraic group $\mathcal{A}_s^0$.  Like any commutative algebraic
group, $\mathcal{A}_s^0$ admits a (unique) \emph{Chevalley decomposition}
\[0 \ra U \oplus T \ra \mathcal{A}_s^0 \ra B \ra 0, \]
where $U$ is a commutative unipotent group (the additive part, even though in positive characteristic
$U$ need not be isomorphic to a power of $\Ga$), $T$ is a linear torus (the multiplicative part)
and $B$ is an abelian variety (the abelian variety part).  
Put $\lambda = \dim U$, $\mu = \dim T$, $\beta = \dim B$. It is standard
to call the reduction \emph{purely additive} if $\lambda = d$.  We shall also say that
the reduction is \emph{atoric} if $\mu = 0$ (this is not standard).  
\\ \\
Potentially good reduction is
atoric \cite[p. 500]{ST} (but no proof is given).  More generally:
\begin{prop}
Let $L/K$ be a finite extension of complete fields, and let $A/K$ be an abelian
variety.  Then $\mu(A/L) \geq \mu(A/K)$ and $\beta(A/L) \geq \beta(A/K)$.  
\end{prop}
\noindent
Proof: Because formation of the N\'eron model commutes with unramified base change, we may
assume without loss of generality that the residue field $k$ of $K$ is separably closed.  Let $S$ be the valuation
ring of $L$ and $\mathcal{A}'/S$ be the N\'eron model of the base change of $A$ to $L$.  On the other
hand, we have $\mathcal{A}/S = \mathcal{A}/R \otimes_R S$, the base change of the N\'eron model of $A/K$
from $R$ to $S$.  If $L/K$ is ramified, these two group schemes need not coincide, but at least
the N\'eron mapping property applied to the smooth $S$-group scheme $\mathcal{A}/S$ and the identity
map on the generic fibres gives us an $S$-morphism $\varphi: \mathcal{A}/S \ra \mathcal{A}'/S$, and in particular
a homomorphism of algebraic groups $\varphi: \mathcal{A}_s/k \ra \mathcal{A}'_s/k$.  The restriction of
$\varphi$ to $T$ must land in the toric part $T'$ of $\mathcal{A}'$, so if the toric rank has 
decreased, there exists some $\Gm \subset \ker(\varphi)$.  This means that
for any $\ell \neq p$, there exists a nontrivial $\ell$-torsion point $P \in \mathcal{A}_s(k)$ mapping to zero in 
$\mathcal{A}'_s(k)$.  Now
$P$ lifts to a unique $\ell$-torsion point $\tilde{P}$ of the generic fibre $A(L)$, and obviously
$0 \neq \varphi(\tilde{P})$; thus $\varphi(\tilde{P})$ reduces to zero in the special fibre of $\mathcal{A}'$, 
but the kernel of the 
reduction map has no $\ell$-torsion (see Proposition 7).  A similar argument
works to show that $\beta$ is non-decreasing, as one checks immediately that if the abelian variety part of
the special fibre over $L$ has smaller dimension, the map $\varphi$ on special fibres kills a nontrivial
abelian subvariety of $\mathcal{A}_s$.
\\ \\
The argument does not, of course, work to show that the additive rank is nondecreasing.  Nevertheless, it
gives the following useful fact:
\begin{prop}
Let $A/K$ be an abelian variety over a complete field and $L/K$ a finite field extension such that
$A/L$ has good reduction.  Then, the preimage of $U(k) \oplus T(k)$ in $A(K)[\tors]$ is mapped
(injectively) into the kernel of reduction $G^1$ of $A/L$.
\end{prop}
\noindent
Finally, we will need the following result, proved in \cite{McCallum}; see also \cite{ELL}.
\begin{thm}(McCallum)
Let $K$ be a complete field with algebraically closed residue field and $A/K$ a $d$-dimensional abelian 
variety with potentially good reduction.  Let $L/K$ be the (unique, minimal, Galois) extension over
which $A$ aquires good reduction. \\
a) The exponent of $\Phi(k)$ divides $[L:K]$. \\
b) If a prime number $\ell$ divides $\#\Phi(k)$, $\ell \leq 2\lambda + 1 \leq 2d+1$.
\end{thm}

\subsection{The proof of parts b) and c)}\ \\
We already have $A(K)[\tors]' \leq [(1+\sqrt{q})^{2d}]$, so it suffices to bound $A(K)[p^{\infty}]$.
By the work of the previous section, there is a filtration on $A(K)$ whose sucessive quotients are $\Phi(k),\  \mathcal{A}^0_s(k)$
and $G_1$.  Intersecting the filtration with $A(K)[p^{\infty}]$, we get subgroups $H_1, \ H_2, \ H_3$ of these
quotients, such that $\#A(K)[p^{\infty}] = \#H_1 \cdot \#H_2 \cdot \#H_3.$  We will discuss each in 
turn. \\ \indent
By Proposition 7, we have $\#H_3 \leq \#(G_1)[\tors] \leq p^{fd[\frac{e}{p-1}]}$. \\ \indent
By Theorem 11, the exponent of the (geometric) component group divides the order of $[L:K^{\unr}]$, where
$L/K^{\unr}$ is the (unique, minimal, Galois) inertial extension over which $A$ acquires good reduction.
By Theorem 6, the Galois group of $L/K^{\unr}$ injects into $\Aut(E[N])$ for any $N \geq 3$ and prime to $p$,
which by definition of $m_p(2d)$ means $\exp(H_1) \leq \exp(\Phi[p^{\infty}]) \ | \ p^{m_p(2d)}$.  But
$A(K)[\tors]$ is generated by $2d$ elements, hence so is $H_1$, hence $\#H_3[p^{\infty}] \ | \ p^{2dm_p(2d)}$.  
\\ \indent
Finally, $\mathcal{A}_s^0$ is an extension of an abelian variety of dimension $d-\lambda$ by a commutative unipotent
group of dimension $\lambda$.  Recall that a commutative unipotent group $U$ over an algebraically field of 
characteristic $p$ is isogenous to a direct product
of finite Witt-vector groups \cite[p. 176]{AGCF} -- in particular, since $\dim U = \lambda \leq d$, we have
$U(\F_q) = U(\F_q)[p^{\lambda}]$ of cardinality $p^{f\lambda}$.  Since $(1+\sqrt{q})^{2a} > q^a$, the worst-case 
scenario for $H_2$ is \emph{good} reduction, and $\#H_2 \leq [(1+\sqrt{q})^{2d}]_p$.     \\ \indent
Taking the product of these three bounds and the bound of part a) gives the proof of part b).
\\ \\
Now assume that $A/K$ has purely additive reduction.  \\ \indent
By Theorem 10b), the only primes dividing $\#\Phi$ are those $\ell \leq 2\lambda +1 \leq 2d+1$.  Since in the purely
additive case $A(K)[\tors]' \hookrightarrow \Phi(k)$, this gives the first assertion of part c).  \\ \indent
As for the second assertion, we want to derive a different bound which is independent of
$p$ and $f$.  Since $p^{2dm_p(2d)} \leq \max \{\#GL_{2d}(\F_3)^{2d}, \#GL_{2d}(\F_5)^{2d}\} \leq
\#GL_{2d}(\F_5)^{2d}$, the bound on $\#H_3$ is already independent of $p$.  The bound on $H_1$ technically
depends on $p$, but in a wonderful way: $\lim_{p \ra \infty} p^{fd[\frac{e}{p-1}]} = 1$, so it may be taken to be
independent of $p$.  Moreover, by Proposition 7a) the exponent of $G_1$ is independent of $f$, hence, reasoning as above,
the order of $H_3$ is independent of $f$.  Finally, the $H_2$ term is handled by Proposition 9, which tells us that $H_2$
injects into ``an $H_3,$'' i.e., into the points of a formal group whose ramification index is bounded by 
$e \cdot \#GL_{2d}(\ell)$, where $\ell$ is any odd prime different from $p$.  This completes the proof of the Main Theorem.
\subsection{Comparisons with the elliptic curve case}\ \\
The case of elliptic curves with purely additive reduction was
studied by Flexor and Oesterl\'e \cite{FO}; they showed $E(K)[\tors] \leq 48e(K/\Q_p)$.  Part c) of the
Main Theorem can be regarded as a generalization of this, and especially of the existence of a bound
independent of $p$ and the residual degree $f$.  
\\ \\
However, our result does not completely generalize the elliptic curve case, because Flexor and Oesterl\'e
do not need to assume that the reduction is potentially good; rather, only that it is \emph{not} split
multiplicative.  The essential fact used is the bound $\#\Phi(\overline{k}) \leq 4$ for the component
group in this case, a well-known consequence of Tate's algorithm.  It is probable that there is an absolute
bound for the exponent (and even the order) of $\Phi$ for an arbitrary atoric abelian variety, but this is only
known for the prime-to-$p$ part.  
\\ \\
Nevertheless, X. Xarles \cite{Xarles} has shown that torsion is strongly bounded among all $d$-dimensional abelian 
varieties over a $p$-adic field whose N\'eron special fibre does not contain a nontrivial split torus.  
Moreover, in the case of purely additive reduction, his bound depends only on $d$ and $e(K/\Q_p)$.
\section{Elliptic Curves over Number Fields}
\subsection{Elliptic curves $E/\Q$} 
\noindent
Let $E/\Q$ be a rational elliptic curve with everywhere
potentially good reduction.  In particular we have potentially good reduction at $2$, so by Theorem 1a) the
odd-order torsion in $E(\Q)$ is a group of order at most $[(1+\sqrt{2})^2] = 5$, so is either
trivial, or $\Z/3\Z$ or $\Z/5\Z$.  But also we have good reduction at $3$, so the prime-to-$3$ torsion
subgroup injects into a group of order at most $[(1+\sqrt{3})^2] = 7$, so is trivial, or is
$\Z/2\Z,\  \Z/2\Z \times \Z/2\Z, \ \Z/4\Z$ or $\Z/5\Z$.  Finally, the prime-to-$5$-torsion injects into 
a group of order at most $[(1+\sqrt{5})^2] = 10$.  Collating this information and comparing with the
statement of Theorem 3, we see that what remains to be shown is that there cannot exist
a rational $5$-torsion point.  \\
\indent
Now we use that $E/\Q$ does not have everywhere good reduction, a well-known theorem of Tate;
it is also an immediate consequence of the elliptic modularity theorem.  Moreover, we claim that there is a 
prime $p \neq 5$ of bad reduction.  Otherwise the (modular!) elliptic curve
would have to have conductor $25$, and since $X_0(25)$ has genus zero, there are no such elliptic curves.  
Thus we get a prime $p$ of bad, so necessarily purely additive reduction, and part c) of the Main Theorem 
applies to show that $E(\Q_p)[5] = 0$, and \emph{a fortiori} $E(\Q)[5] = 0$.
\\ \indent
That all these groups occur already among elliptic curves $E/\Q$ with $j$-invariant $0$
or $1728$ is a very classical fact, which can readily be reestablished by looking at the 
Weierstrass equations $y^2 = x^3 +Dx$ ($j= 0$) and $y^2 = x^3 + D$ ($j = 1728$).
This completes the proof of Theorem 3.
\subsection{A miscellany of bounds}
Note that the bounds of Corollary 2 are exponential in both the dimension $d$ and the number field degree
$n$.  That the bound is exponential in $d$ cannot be helped: if $E/\Q$ is a CM elliptic curve with
$E(\Q)[\tors] \cong \Z/6\Z$, then $E^d$ is a $d$-dimensional abelian variety with integral moduli 
and a torsion group of size $6^d$.  Since $6^d$ is, with respect to $d$, the leading term in
the bound of Corollary 2a), this bound looks quite good when $d$ is large compared to $n$.
\\ \\
On the other hand, there is no reason to believe that the optimal bound for the torsion subgroup
of an elliptic curve $E$ over a degree $n$ number field should be exponential in $n$, even if we
do not restrict to integral moduli.  Indeed, there is evidence to the contrary: Hindry and
Silverman \cite{HS} showed that by restricting to elliptic curves with everywhere \emph{good} reduction over
a number field of degree $d$, the torsion subgroup is at most $O(d\log d)$ (with an explicit constant).
Moreover, Silverberg \cite{Silverberg 1992} gives bounds for CM abelian varieties of every dimension: when $d$ is large
compared to $n$, our Corollary 2 still gives better bounds,\footnote{The bounds of Corollary 2 restricted to CM varieties
also supercede those of \cite{Van Mulbregt}.} but when $n$ is large compared to $d$, Silverberg's bounds
are much, much better.  
\\ \\
In the case of elliptic curves an  independent, and more elementary, proof of Silverberg's theorem was given by 
Prasad and Yogananda 
\cite{PY}.  Here is the bound: 
\begin{thm}(Silverberg, Prasad-Yogananda)
Let $E/K$ be an $\OO$-CM elliptic curve over a number field $K$ of degree $n$.  Then
\[\varphi(e) \leq \delta(\#\mu(\OO))n, \]
where $e$ is the exponent of $E(K)[\tors]$, and $\delta$ is $\frac{1}{2}$ if $K$ contains the CM field and is 
$1$ otherwise.  Moreover, if $K$ does \emph{not} contain the CM field, then the same bound holds for the order of the 
torsion subgroup.
\end{thm}
\noindent
For example, when $n =5$, Corollary 2b) says that the largest prime dividing the torsion subgroup is at most
$43$, whereas Theorem 11 tells us that the largest prime divisor is at most $31$.  
\\ \\
The values of $N$ for which a CM elliptic curve $E/K$ can have full $N$-torsion are especially
restricted.  Indeed we state without proof the following result (it is a consequence of the main 
theorem of complex multiplication):
\begin{prop}
Let $E/K$ be a CM elliptic curve over a number field.  If for some positive integer $N \geq 3$, 
$(\Z/N\Z)^2 \hookrightarrow E(K)[\tors]$,
then $H(N) \subset K$, where $H(N)$ is the ring class field of conductor $N$.
\end{prop}
\noindent
Remark: Since $\Q(\zeta_N)$ is contained in $H(N)$, it follows that $K$ must contain the $N$th roots of unity; of course this 
latter
containment holds for \emph{all} elliptic curves by the Galois equivariance of Weil's $e_N$-pairing.  The proposition tells us 
that in the CM case, for $N \geq 3$, $\Q(E[N])$ will typically be much larger than $\Q(\zeta_N)$. 
\\ \\
A moment's reflection reveals that it ought to be the case that the elliptic curves defined over a degree $n$ number
field with the largest torsion subgroups are those which arise via base change from a smaller field, or especially, from $\Q$.
A beautifully precise confirmation of this for CM elliptic curves is given by the following result of J.L. Parish \cite{Parish}.
\begin{thm}(Parish)
Let $E/K$ be a CM elliptic curve which is \emph{minimally} defined over $K$ -- i.e.,  $K = \Q(j(E))$ -- and let 
$F$ be the compositum
of $K$ with the CM field.  Then the order of a torsion point on $E(F)$ is $1, \ 2, \ 3, \ 4,$ or $6$.  Moreover, the torsion
subgroup of $E(K)$, if not cyclic, is of the form $\Z/2\Z \times \Z/2\Z$.
\end{thm}
\noindent
\subsection{An algorithm for enumerating torsion subgroups of CM elliptic curves over number fields of prime degree}
\noindent \ \\
Suppose $n$ is a prime number.  We give a procedure for enumerating the torsion subgroups of CM elliptic curves
over degree $n$ number fields.  
\\ \\
We will assume first that $n$ is odd; the case $n=2$ will be settled at the end of the section.
Because $n$ is odd, Proposition 11 (or even the remark following it) allows us to conclude that $E(K)[\tors]$ is of the form
$\Z/N\Z \times \Z/a\Z$, where $a \ | \ gcd(N,2)$.  In particular, the odd-order torsion subgroup is cyclic.
\\ \\
First of all, the groups $\Z/n\Z$ for $n \in \{1,2,3,4,6\}$ and $\Z/2\Z \times \Z/2\Z$ occur over $\Q$, hence
will by base change occur for all degrees $n$.  (It follows easily from the strong boundedness of torsion on elliptic
curves that if $E/K$ is an elliptic curve over a number field and $n$ a positive integer, then there are at most finitely many 
degree $n$ field extensions $L/K$ such that $E(K)[\tors] \subsetneq E(L)[\tors]$.)  
\\ \\
Now the key point is that since $n$ is prime, either $\Q(j(E)) = \Q$ or $\Q(j(E)) = K$.  But Theorem 12 tells us precisely
that in the latter case we will not get any torsion subgroups other than those of the last paragraph.  Thus, it suffices
to restrict attention to the $13$ $\Z$-integral CM $j$-invariants, namely $0, \ 1728, \ 2^43^35^3, \ -2^{15}35^3, \ 2^33^311^3, \\ \ -3^35^3, \ 3^35^317^3, \ 2^65^3, \ -2^{15}, \ 
-2^{15}3, \ -2^{18} 3^3 5^3, \ -2^{15} 3^35^311^3, \ -2^{18}3^3 5^3 23^3 29^3.$
\\ \\
The first step of the algorithm is to apply Corollary 2 and/or Theorem 11 to compile a short list of candidate values of
$N$.  Then, for each such value of $N$ not in the set $\{1,2,3,4,6\}$ and each $j_i$, $1 \leq i \leq 13$, there is a very direct
way to decide whether or not there exists an elliptic curve $E$ over a degree $n$ field with $j(E) = j_i$ and
$\Z/N\Z \subset E(K)[\tors]$: namely, we compute the zero-dimensional subscheme $T_{N,j_i}$ of $X_1(N)/\Q$ obtained
as the preimage of the point $j_i \in \mathbb{A}^1(\Q)$ under the map
\[X_1(N) \ra X(1) \stackrel{j}{\ra} \PP^1. \]
As we will shortly describe, $T_{N,j_i}$ may be represented as a zero-dimensional ideal in a multivariable
polynomial ring over $\Q$.  Its radical thus determines a finite \'etale $\Q$-algebra -- i.e., a product of number fields 
$\prod_{i=1}^I K_i$ -- and we define the \emph{degree sequence} $D(N,j_i)$ to be the $I$-tuple giving the degrees over
$\Q$ of each such number field $K_i$.  Thus, there exists an $N$-torsion point on a CM elliptic curve
with $j$-invariant $j_i$ defined over a number field of degree $n$ if and only if $n$ appears in the degree sequence
$D(N,j_i)$.  Moreover, if it does, this procedure allows us to actually write down the corresponding elliptic curve,
and its full torsion subgroup may be determined.  
\\ \\
We now explain how $T_{N,j_i}$ may be computed.  Suppose we have $(E,P)/K$, where $K$ is a field of characteristic zero, $E/K$ is 
an elliptic curve and $P \in E(K)$ is a
point of exact order $N \geq 4$.  Starting with any given Weierstrass model for $E$, one can make a linear change of variables
to get the \emph{Kubert normal form}
\[E(b,c): y^2 + (1-c)xy -by = x^3-bx^2, \ b, \ c \in K, \]
with $j$-invariant
\[j(b,c) =   \frac{\left(((1-c)^2-4b)^2+24b(1-c)\right)^3}{b^3(((1-c)^2-4b)^2+8(1-c)^3-27b-9(1-c)((1-c)^2-4b))}. \]
The function $j(b,c)$ is the composite function we want, since $\Q(X_1(N)) = \Q(b,c)$.  As for the relations
between $b$ and $c$, they can be computed directly from the group law, i.e., by writing down the equations corresponding
to $[N]P = 0, \ [i]P \neq 0$ for $0 < i < N$.  When $X_1(N)$ has genus zero, the $j$-function on $X(1)$ can be given
as a rational function in some parameter $T$ on $X_1(N)$, so that to compute $T_{N,j_i}$ amounts to factoring a one-variable
polynomial over $\Q$.  Indeed, the parameterizations $T$ are ``classical,'' having been computed by Kubert \cite{Kubert}
in his exploration of Ogg's conjecture -- i.e., Mazur's theorem -- on the torsion subgroup of $E/\Q$.  Work of
Reichert \cite{Reichert} demonstrates that computing the relations between $b$ and $c$ is not essentially more difficult in the 
higher genus situation; indeed he computes such relations (his ``raw forms'') when $N = 11, \ 13, \ 14, 15, 16, \ 17, 18$;
note that $X_1(17)$ has genus $5$.  Since $[(1+2^{3/2})^2] = 14$, Corollary $2$ tells us that the largest prime value of
$N$ that we need consider in order to prove Theorem 4 is $13$.  Thus, in order to prove Theorem 4 for $n = 3$, it suffices to 
type the equations
of \cite{Kubert} and \cite{Reichert} and the above list of $j$-invariants into a computer algebra system and record the
degree sequences.  They are as follows: (the $j$-invariants are listed in the same order as above; to save space,
when a degree sequence $(d_1,\ldots,d_I)$ is repeated $a > 1$ times, we denote it by $(d_1,\ldots,d_I)^a$):
\\ \\
$\Z/4\Z$: $(2), (1,2), (2,4), (6), (1,1,4), (2,2,2), (2,4)^2, (6)^5$ \\
$\Z/5\Z$: $(4), (2,4), (12)^2, (4,8), (12)^3, (4,8)^2, (12)^3$ \\
$\Z/6\Z$: $(1,3), (2,4), (1,2,3,6), (3,9), (4,8)^3, (2,2,4,4), (6,6), (12)^4$ \\
$\Z/7\Z$: $(2,6), (12), (6,18)^2, (24), (3,21)^2, (24)^2, (6,18), (24)^3$ \\
$\Z/8\Z$: $(8), (4,8), (8,16), (24), (4,4,16), (4,4,8,8), (4,4,16), (8,16), 
(24)^5$ \\
$\Z/9\Z$: $(3,9), (18), (9, 27), (36)^4, (6,12,18)^2,
(36)^4$ \\
$\Z/10\Z$: $(12), (2,4,4,8), (12,24), (36), (4,8,8,16), (12,24)^5, (36)^3$ \\
$\Z/11\Z$: $(20), (30), (60)^3, (10,50)^3, (5,55), (10,50)^2, (60)^2$ \\
$\Z/12\Z$: $(4,8), (8,16), (4,8,12,24), (48), (8,8,32), (16, 16, 16), (16, 32), 
(8,8,16,16), (24, 24), (48)^4$
$\Z/13\Z$: $(4,24), (6,36), (12,72), (84), (12,72), (84)^5,
(12,72), (84)^2$
\\ \\
Kubert also gives defining equations for genus zero torsion structures
of the form $\Z/2\Z \times \Z/2M\Z$, so we give those degree sequences as well:
\\ \\
$\Z/2\Z \times \Z/4\Z$: $(4), (2,4), (4,8), (12), (2,2,8), (4,4,4), (4,8)^2, 
 (12)^5$ \\
$\Z/2\Z \times \Z/6\Z$: $(2,6), (4,4,4), (2,2,2,6,6,6), (24), (8,8,8)^3, 
(4,4,4,4,4,4), (6,6,12), (24)^4$ \\
$\Z/2\Z \times \Z/8\Z$: $(16), (8,8,8), (8,8,16,16), (48), (4,4,8,16,16),
(4,4,4,4,8,8,16), (8,8,16,16)^2, (48)^5$
\\ \\
Remark: For fixed $N$ and for all $j_i$ except $0$ and $1728$,
the sums of the degree sequences are constant, and equal to the degree of
the natural map $\varphi: X_1(N) \ra X_1(1)$, whereas the total degrees at $j = 0$
and $j = 1728$ are always smaller: this exhibits the well-known fact that
$\varphi$ is ramified only at $j = 0$ and $j = 1728$ (and possibly at the cusps).
\\ \\
To complete the proof of Theorem 4 for $n = 3$ we need only look at the three instances of ``3''
among the degree sequences and compute the full torsion subgroup of each of the
corresponding elliptic curves.  These three elliptic curves are
given, in Kubert normal form, by the following values of $b$ and $c$:
\[K = \Q[d]/(d^3-d^2-2d+1), \ b = 2d-1, \ c = d^2-d, \ E(K)[\tors] \cong \Z/14\Z. \]
\[K = \Q[d]/(d^3-15d^2+12d+1),\  b = 14d^2-12d-1, \ c = d^2-d, \ 
E(K)[\tors] \cong \Z/14\Z. \]
\[K = \Q[f]/(f^3-3f^2+1), \ b = 13f^2-f-5, \ c = 2f^2-1, \ E(K)[\tors] \cong \Z/9\Z. \]
For the case $n=2$, we must contemplate the possibility that $(\Z/N\Z)^2 \subset E(K)[\tors]$
for $N = 3,\ 4$ or $6$.  
\\ \\
For $N =3$, the modular curve $X(3)/\Q(\zeta_3)$ is again ``parameterizable'' (i.e., of genus zero): 
$\Q(\zeta_3)(X(3)) \cong
\Q(\zeta(3))(\lambda)$, and the equation of the universal elliptic curve over the generic
fibre is especially well-known: it is the \emph{Hesse curve}
\[E_{\lambda}: X^3+Y^3+Z^3 + \lambda XYZ = 0. \]
Here, if we take $O = [1:-1:0]$, then a basis for $E[3]$ is given by
$S = [1,\zeta_3,0], \ T = [1,0,-1]$.  The $j$-invariant of the Hesse curve is
\[J_{3 \times 3}(\lambda) = 
\frac{\lambda^{12}-648\lambda^9+139968\lambda^6-10077696\lambda^3}{-\lambda^9-81\lambda^6-2187\lambda^3-19683}. \]
Using this, one sees that the unique CM elliptic curve over $\Q(\zeta_3)$ with full $3$-torsion is
the Fermat curve $X^3+Y^3+Z^3=0$, whose full torsion group is $(\Z/3\Z)^2$.  In particular there is no CM elliptic
curve over a quadratic field with full $6$-torsion.\footnote{Ironically, since $X(6)$ has genus one and an order six
automorphism with a fixed point, it too must be the Fermat curve $X^3+Y^3+Z^3 = 0$.}
\\ \\
The case of $(\Z/4\Z)^2 \subset E(K)[\tors]$ can be ruled out by analyzing the three
instances of ``$2$'' in the degree sequences for $(\Z/2\Z \times \Z/4\Z)$ (or indeed by
looking at $\Gamma(4)$ level structure, which has an equally well-known parameterization).  For the sake of
variety, we give the following alternate argument: if $E/K$ is any elliptic curve defined
over a quadratic number field with full $4$-torsion, then necessarily $K = \Q(\zeta_4)$.
So $5$ splits in $K$, and if $E$ has integral moduli, the $4$-torsion would have to reduce
injectively into an elliptic curve over $\F_5$, which however has order at most
$[(1+\sqrt{5})^2] = 10$.  Thus it is impossible for an IM elliptic curve over a quadratic
field to have full $4$-torsion.
\\ \\
This completes the proof of Theorem 4.
\section{Abelian surfaces over $\Q$}
\subsection{Bounds obtained from the Main Theorem} 
Let $A/\Q$ be an IM abelian surface.  As usual, applying part a) of the Main Theorem 
prime by prime
leads to a short list of possible orders of torsion groups.  Indeed, in this case, the odd order torsion injects 
into a group of order at most $\#A(\F_2) \leq [(1+\sqrt{2})^4] = 33$
and the prime-to-$3$ torsion injects into a group of order atmost $\#A(\F_3) \leq [(1+\sqrt{3})^4] = 55$.
This implies that the the possible orders of torsion groups are of the form $2^a \cdot y$, where
$0 \leq a \leq 5$ and $y$ lies in the set
\newcommand{\odd}{\operatorname{odd}}
\begin{equation}
1,3,5,7,3^2,11,13,3\cdot5,17,19,3\cdot 7,23,5^2,3^3,29,31,3\cdot 11. 
\end{equation}
\subsection{Some Tate-Honda Theory}
Suppose one wishes to enumerate all the possible values of $\#A(\F_q)$, where $A/\F_q$
is a $d$-dimensional abelian variety.  Recall that for two abelian varieties
$A_1, \ A_2$ over $\F_q$, the following are equivalent: a) that they are isogenous,
b) that their Frobenius characteristic polynomials coincide, c) that they have the
same number of rational points over every finite field extension.  Thus, to perform
the enumeration, it is enough to know the set of all Frobenius polynomials
$P(T)$ of $d$-dimensional $A/\F_q$: just evaluate at $T=1$.
\\ \\
This problem -- namely, which polynomials arise as Frobenius polynomials? -- is 
addressed by the theory of Tate and Honda.  The definitive introduction to this
theory is still to be found in Waterhouse's thesis \cite{Waterhouse}; here we
will just give an ``explicit formula'' for $\#A(\F_p)$ where $A$ is an abelian surface.
\begin{prop}
There are three ``types'' of abelian surface $A/\F_p$, whose
Frobenius polynomials are as follows: \\
a) Type I: $P_A(T) = (T^2-a_1T+p)(T^2-a_2T+p)$,
where $a_1,\ a_2$ are integers such that $|a_1|, \ |a_2| < 2\sqrt{p}$;
$\#A(\F_p) = (p+1-a_1)(p+1-a_2)$. \\
b) Type II: $P_A(T) = (T^2-p)^2$;
$\#A(\F_p) = (p-1)^2$. \\
c) Type III: $P_A(T) = T^4-2aT^3+(a^2+2p-2db^2)T^2-2apT + p^2$,
where $d > 1$ is a squarefree integer, $a, \ b \in \frac{1}{2}\Z$ are such that
$a+b\sqrt{d}$ is in the ring of integers of $\Q(\sqrt{d})$, $b \neq 0$
and $|a|+|b|\sqrt{d} < 2\sqrt{p}$; 
$\#A(\F_p) = (p+1)^2+(a-1)^2-2ap-db^2-1$.
\end{prop}
\noindent
Proof: If $A/\F_q$ is an abelian surface, then either
$A \sim_{\F_q} E_1 \times E_2$ or $A$ is $\F_q$-simple.
The former case is Type I, so it suffices to consider the $\F_q$-simple (henceforth called
``simple,'' which is not to be confused with geometrically simple) case.
\\ \indent
If $A/\F_q$ is a simple abelian surface, its Frobenius polynomial
$P_A(T)$ is a quartic Weil $q$-polynomial, i.e., a polynomial in $\Z[T]$ 
whose roots have norm $\sqrt{q}$ under every Archimedean valuation.  Moreover,
$P_A$ is either irreducible or is of the form $Q(T)^2$ for $Q(T)$ an irreducible
quadratic.  Over a general finite field $\F_q$ it is somewhat intricate to describe which
Weil polynomials of the second type correspond to abelian surfaces, but over $\F_p$
this can only happen if the field $\Q(\pi)$ generated by a Frobenius root $\pi$ 
is real, i.e., $\pi = (\pm)\sqrt{p}$.  This is Type II.  Otherwise
$\Q(\pi)$ is a quartic CM field which is best understood in terms of the real
quadratic subfield $\Q(\beta)$, where $\beta = \pi + \frac{p}{\pi}$.  Indeed, the
condition on $\beta$ that it be ``the $\beta$'' of some quartic Weil $p$-number
$\pi$ is just that it be an irrational real quadratic integer $\beta = a+b\sqrt{d}$
which has norm strictly less than $2\sqrt{p}$ in both Archimedean valuations, i.e., 
$|a|+|b|\sqrt{d} < 2\sqrt{p}$.  The corresponding
$\pi$ is then a solution of $T^2-\beta T+p = 0$.  This is Type III.
\subsection{The proof of Theorem 5}
Using Proposition 14, we record $\#A(\F_p)$ for $p = 2, \ 3, \ 5$:
\newcommand{\OR}{\operatorname{or}}
\begin{fact}
Let $A/\F_p$ be an abelian surface over the finite field $\F_p$, $p \leq 5$.  Then:
\[\#A(\F_2) = 1-16, \ 19, \ 20, \ \OR \  25. \]
\[\#A(\F_3) = 1-16, \ 18-25, \ 28-30, \ 34-36, \ 42, \ \OR  \ 49. \]
\[\#A(\F_5) = 4, \ 6-50, \ 52-56, \ 58-64, \ 69-72, \ 79-81, \ 90, \ \OR  \ 100. \]
\end{fact}
\noindent
 Note the
sparsity for $p=2$ and $p=3$: by the time we get to $p=5$ there is such an enormous interval
of assumed values that nothing further is ruled out.\footnote{In an earlier draft I wrote ``lucky sparsity,'' but \cite{DH} gives evidence that the smaller size of the 
\emph{central interval}
$I_{d,q}$ -- i.e., the largest symmetric interval $I$ centered at $q^{d}+1$ such that if $N \in I$, then there exists
$A/\F_q$ of dimension $d$ with $\#A(\F_q) = N$ -- for $q = 2$ and $3$ than for other prime powers is a general phenomenon.}
\\ \\
Working prime by prime, we now use Fact 15 to eliminate many of the values 
from the list $2^a \cdot y$, $1 \leq a \leq 5$, $y$ in the list (6) (e.g., we can in
fact take $a \leq 4$), and in doing so we arrive at the list (4) of Theorem 5.
\subsection{Some attained values of $N$}
\noindent 
Let us discuss which values of the list (4) we know do arise.
\\ \\
If $E_1, \ E_2 / \Q$ are two IM elliptic curves, then $A := E_1 \times E_2/\Q$ is an IM
abelian surface, and $A(\Q)[\tors] = E_1(\Q)[\tors] \times E_2(\Q)[\tors]$.  Thus, using Theorem 3,
the orders of $A(\Q)[\tors]$ arising in this way are: $1-4, \ 6, \ 8, \ 9, \ 12 \ 18, \ 24, \ 36$.
\\ \\
\newcommand{\Res}{\operatorname{Res}}
If $E/K$ is an IM elliptic curve over a quadratic field, then $A: = \Res_{K/\Q}(E)$, the Weil 
restriction of $E$ from $K$ to $\Q$, is an IM abelian surface with $A(\Q)[\tors] = E(K)[\tors]$.  Thus the classification
of torsion subgroups on IM elliptic curves over quadratic fields of \cite{Z1} becomes relevant: it
turns out that, in addition to the groups of Theorem 4a), there arise $\Z/5\Z$ and $\Z/8\Z$.  So 
in this way, we get additional orders $5, \ 7, \ 10$  for $A(\Q)[\tors]$.
\\ \\
The abelian surface $J_1(13)/\Q$ -- i.e., the Jacobian of the modular curve $X_1(13)$ -- has integral
moduli and $J_1(13)(\Q)[\tors] \cong \Z/19\Z$.  Thus $19$ arises.
\\ \\
The abelian surface $J_1(16)/\Q$ has integral moduli and $J_1(16)(\Q)[\tors] \cong \Z/2 \times
\Z/10$.  Thus $20$ arises.
\\ \\
So to complete the classification of orders of torsion subgroups of IM abelian surfaces $A/\Q$,
it remains to decide which of the following eleven values occur:
\begin{equation}
11, \ 13, \ 14, \ 15, \ 22, \ 25, \ 28, \ 30, \ 48, \ 60, \ 72.
\end{equation}

\end{document}